\documentclass[11pt]{amsart}

\usepackage{amsmath,amssymb,amsthm}
\usepackage{hyperref}

\numberwithin{equation}{section}

\newtheorem{prop}{Proposition}[section]

\newtheorem{lemm}[prop]{Lemma}

\newtheorem{rema}[prop]{Remark}

\providecommand{\func}[1]{\operatorname{#1}}

\title[A Concavity Inequality for Hessian Quotient Equations]
{A Concavity Inequality for Hessian Quotient Equations}

\author{Yi-Lin Tsai}
\address{Department of Mathematics and Statistics, McMaster University, Hamilton, Ontario, Canada}
\email{tsaiy11@mcmaster.ca}

\begin{document}

\begin{abstract}
We prove a concavity inequality for the Hessian quotient operator
$\sigma_k/\sigma_{k-1}$ using a change of basis for symmetric polynomials.
As a consequence, the interior $C^2$ estimate for convex solutions of the
$\sigma_k/\sigma_{k-1}$ equation holds without the additional structural
concavity assumption imposed in previous work. The method also yields
several useful inequalities for elementary symmetric polynomials.
\end{abstract}

\maketitle

\section{Introduction}

Concavity inequalities for elementary symmetric functions play an important
role in second-order estimates for Hessian and curvature equations.
Guan--Ren--Wang \cite{GRW} used such an inequality implicitly in their global
$C^{2}$ estimate for convex solutions of curvature equations. Chu
\cite[Lemma~3.3]{Chu} later stated it explicitly in a form close to the one
here. Related concavity inequalities were developed by Ren and Wang for the
$(n-1)$- and $(n-2)$-Hessian equations on the corresponding G\aa rding cones;
see \cite{RW,RW2,RWNotes}. Lu \cite{LuHyperbolic} obtained a related
inequality for semi-convex solutions of Hessian equations. More recently,
Zhang \cite[Lemma~1.1]{Zhang} proved a closely related inequality for
$k$-Hessian equations under a semi-convexity assumption.

For Hessian quotient equations, a stronger concavity inequality with an
additional positive term is needed. Guan and Sroka \cite{GS} proved such an
inequality for positive Hessian quotient operators. This inequality gives
the key Jacobi inequality used in Lu's interior $C^{2}$ estimate for Hessian
quotient equations \cite{Lu24}. In particular, Lu studied
\begin{equation*}
\frac{\sigma_n}{\sigma_l}(D^{2}u)=f
\end{equation*}
and obtained interior $C^{2}$ estimates for $l=n-1,n-2$.

In our previous work \cite{LuTsai}, we considered the general Hessian
quotient equation
\begin{equation*}
\frac{\sigma_k}{\sigma_l}(D^{2}u)=f,\qquad l=k-1,\ k-2,
\end{equation*}
and proved an interior $C^{2}$ estimate for convex solutions under an
additional structural concavity assumption. The main result of this paper
establishes this concavity inequality when $l=k-1$. Consequently, the
interior $C^{2}$ estimate for the $\sigma_k/\sigma_{k-1}$ equation holds
without this additional assumption.

Let
\begin{equation*}
F=\frac{\sigma_k}{\sigma_{k-1}}(\lambda),\qquad
\lambda_1\geq\lambda_2\geq\cdots\geq\lambda_n>0,\qquad 2\leq k\leq n-1.
\end{equation*}
We prove the following.

\begin{lemm}
\label{mainlemma} There exists $\varepsilon>0$ depending only on $n$ and $k$, and
$M=M(\varepsilon,n,k,F)$, such that if $\lambda_1>M$, then, for any
$\xi\in\mathbb{R}^n$,
\begin{equation*}
\frac{1}{F}\left(\sum_i F^{ii}\xi_i\right)^2 -\sum_{i,j}F^{ii,jj}\xi_i\xi_j
+2\sum_{i>1}\frac{F^{ii}\xi_i^2}{\lambda_1} -\frac{F^{11}\xi_1^2}{\lambda_1}
\geq \varepsilon\frac{F^{11}\xi_1^2}{\lambda_1},
\end{equation*}
where
\begin{equation*}
F^{ii}=\frac{\partial F}{\partial\lambda_i},\qquad F^{ii,jj}=\frac{%
\partial^2F}{\partial\lambda_i\partial\lambda_j}.
\end{equation*}
\end{lemm}

Thus, the structural assumption in \cite{LuTsai} holds for $l=k-1$. Together
with \cite{LuTsai}, Lemma~\ref{mainlemma} gives the interior $C^2$ estimate
for convex solutions of the $\sigma_k/\sigma_{k-1}$ equation without this
additional assumption.

\begin{rema}
After this work was completed, Li and Wu \cite{LiWu} independently obtained
a related concavity inequality under different conditions. Their different
approach also applies to $\sigma_k/\sigma_{k-2}$.
\end{rema}

We now give a brief outline of the proof. A determinant computation reduces
Lemma~\ref{mainlemma} to a scalar inequality involving elementary symmetric
polynomials. The endpoint cases $k=2$ and $k=n-1$ are proved directly. For
$3\leq k\leq n-2$, we further reduce the problem to the positivity of a cubic
polynomial
\begin{equation*}
N(x)=A_3x^3+A_2x^2+A_1x+A_0.
\end{equation*}
Using Rolle's theorem and reciprocal variables, we reduce the terms
involving $e_{k-1},\ldots,e_{k-4}$ to the first four elementary symmetric
polynomials $e_1,\ldots,e_4$.

The main idea is to \emph{change the basis} of symmetric polynomials.
This gives useful inequalities for elementary symmetric polynomials that
are not transparent from the Newton--Maclaurin inequalities alone.
The coefficients $A_i$ also have a multi-affine structure when expressed
in terms of certain ratios of elementary symmetric polynomials.
Together, these facts reduce the positivity problem to checking values at
the vertices of suitable rectangles, giving a direct proof of the desired
concavity inequality.

We also give an independent computational verification using SageMath
\cite{S}. Expanding the coefficients $A_i$ in the monomial basis reduces
their positivity to coefficient positivity together with a few elementary
comparisons.

The paper is organized as follows. Section~2 gives some preliminary
formulas. Section~3 proves the symmetric polynomial inequalities used later.
Section~4 reduces the main estimate to a scalar inequality, and Section~5
treats the endpoint cases. Section~6 reduces the remaining cases to a cubic
polynomial. Section~7 gives an independent computational verification using the
monomial basis. Finally, Section~8 gives a direct proof using the
multi-affine structure and the inequalities from Section~3.

\section{Preliminaries}

We collect some basic properties of the Hessian and Hessian quotient
operators.

For $\lambda=(\lambda_1,\ldots,\lambda_n)\in \mathbb{R}^n$, we denote
\[
(\lambda|i)=(\lambda_1,\ldots,\lambda_{i-1},\lambda_{i+1},\ldots,\lambda_n)
\in \mathbb{R}^{n-1}.
\]
Thus, $(\lambda|i)$ is obtained from $\lambda$ by deleting its $i$-th
component. Similarly, $(\lambda|ij)$ is obtained by deleting its $i$-th and
$j$-th components.

We recall some basic formulas for the elementary symmetric functions; see,
for instance, \cite{GRW, LT}.

\begin{lemm}
\label{Sigma_k-Lemma-0} For any $\lambda =(\lambda _{1},\cdots ,\lambda
_{n})\in \mathbb{R}^{n}$, we have 
\begin{align*}
\sigma_k(\lambda) &=\lambda_i\sigma_{k-1}(\lambda|i)+\sigma_k(\lambda|i), \\
\sum_i\lambda_i\sigma_{k-1}(\lambda|i) &=k\sigma_k(\lambda), \qquad
\sum_i\sigma_k(\lambda|i)=(n-k)\sigma_k(\lambda), \\
\sigma_{k-1}^2(\lambda|ij) -\sigma_k(\lambda|ij)\sigma_{k-2}(\lambda|ij)
&=\sigma_{k-1}(\lambda|i)\sigma_{k-1}(\lambda|j)
-\sigma_k(\lambda)\sigma_{k-2}(\lambda|ij).
\end{align*}
\end{lemm}

When there is no ambiguity, we write
\[
\sigma_{k;i}=\sigma_k(\lambda|i),\qquad
\sigma_{k;ij}=\sigma_k(\lambda|ij).
\]
We also recall Newton's inequality
\begin{equation*}
\sigma _{k-1}^{2}\geq \frac{k(n-k+2)}{(k-1)(n-k+1)}\sigma _{k}\sigma _{k-2}.
\end{equation*}

\medskip

Define the G{\aa}rding cone
\[
\Gamma_k=\{\lambda\in \mathbb{R}^n:\sigma_j(\lambda)>0,\ 1\leq j\leq k\}.
\]

If $F$ is a symmetric function of the eigenvalues
$\lambda=(\lambda_1,\ldots,\lambda_n)$ of $D^2u$, we also regard $F$ as
a function of $D^2u$ and write
\[
F^{pq}=\frac{\partial F}{\partial u_{pq}},\qquad
F^{pq,rs}=\frac{\partial^2F}{\partial u_{pq}\partial u_{rs}}.
\]
We also recall the following formulas for the Hessian quotient operator;
see \cite{Lu24}.

\begin{lemm}
\label{Fijkl} Let $1\leq k<\ell \leq n$, let $F=\frac{\sigma _{\ell }}{%
\sigma _{k}}$ and let $u$ be a convex function. Suppose $D^{2}u$ is
diagonalized at $x_{0}$. Then at $x_{0}$, we have 
\begin{equation*}
F^{pp}=\partial _{\lambda _{p}}\left( \frac{\sigma _{\ell }}{\sigma _{k}}%
\right) =\frac{\sigma _{\ell -1;p}}{\sigma _{k}}-\frac{\sigma _{\ell }\sigma
_{k-1;p}}{\sigma _{k}^{2}},
\end{equation*}%
\begin{equation*}
F^{pp,qq}=\partial _{\lambda _{q}}\partial _{\lambda _{p}}\left( \frac{%
\sigma _{\ell }}{\sigma _{k}}\right) 
\begin{cases}
-2\frac{\sigma _{\ell -1;p}\sigma _{k-1;p}}{\sigma _{k}^{2}}+2\frac{\sigma
_{\ell }\left( \sigma _{k-1;p}\right) ^{2}}{\sigma _{k}^{3}},\quad & p=q. \\ 
\frac{\sigma _{\ell -2;pq}}{\sigma _{k}}-\frac{\sigma _{\ell -1;p}\sigma
_{k-1;q}}{\sigma _{k}^{2}}-\frac{\sigma _{\ell -1;q}\sigma _{k-1;p}}{\sigma
_{k}^{2}}-\frac{\sigma _{\ell }\sigma _{k-2;pq}}{\sigma _{k}^{2}}+2\frac{%
\sigma _{\ell }\sigma _{k-1;p}\sigma _{k-1;q}}{\sigma _{k}^{3}}, & p\neq q,%
\end{cases}%
\end{equation*}
\end{lemm}

Finally, recall that
\[
F=\frac{\sigma_k}{\sigma_{k-1}}
\]
is concave in $\Gamma_k$. In particular, when all eigenvalues are positive,
we have
\begin{equation}
C_{1}\left( n,k,F\right) <\lambda _{k}<C_{2}\left( n,k,F\right) .
\label{kbound}
\end{equation}

\section{Symmetric polynomial inequalities}

In this section, we establish several inequalities for $P,Q,R$ that will be
used later in the proof of the main estimate. The key inequalities are
obtained by changing between different bases of symmetric polynomials.

Let $e_j$ denote the $j$-th elementary symmetric polynomial in $k$ positive
variables. Define 
\begin{equation*}
P=\frac{e_1e_3}{e_2^2},\qquad Q=\frac{e_2}{e_1^2},\qquad R=\frac{e_2e_4}{%
e_3^2}.
\end{equation*}
We use the standard change-of-basis formulas for symmetric functions; see 
\cite[Chapter I]{M}. The computations below were also verified using
SageMath~\cite{S}.

\subsection{Inequalities from power sums}

Let $p_{j}=\sum_{i}x_{i}^{j}$ be the $j$-th power sum. By Cauchy--Schwarz, 
\begin{equation}
p_{2}^{2}\leq p_{1}p_{3}.  \label{Power1}
\end{equation}%
Using the following expansions, 
\begin{align*}
p_{1}& =e_{1},\qquad p_{2}=e_{1}^{2}-2e_{2},\qquad
p_{3}=e_{1}^{3}-3e_{1}e_{2}+3e_{3}, \\
p_{4}& =e_{1}^{4}-4e_{1}^{2}e_{2}+4e_{1}e_{3}+2e_{2}^{2}-4e_{4},
\end{align*}%
we obtain from \eqref{Power1} 
\begin{equation}
e_{1}^{2}e_{2}+3e_{1}e_{3}-4e_{2}^{2}\geq 0.  \label{Gap0}
\end{equation}%
Equivalently, 
\begin{equation*}
1\geq (4-3P)Q.
\end{equation*}%
We also use the inequality
\begin{equation*}
p_{4}\leq p_{2}^{2}.
\end{equation*}%
After substitution and simplification, this gives 
\begin{equation*}
2P^{2}QR\geq 2P-1.
\end{equation*}

\subsection{Inequalities from Schur polynomials}

We next use the positivity of Schur polynomials. We use the expansions 
\begin{equation*}
s_{(2,2)}=%
\begin{vmatrix}
e_{2} & e_{3} \\ 
e_{1} & e_{2}%
\end{vmatrix}%
\geq 0,\qquad s_{(3,3)}=%
\begin{vmatrix}
e_{2} & e_{3} & e_{4} \\ 
e_{1} & e_{2} & e_{3} \\ 
e_{0} & e_{1} & e_{2}%
\end{vmatrix}%
\geq 0.
\end{equation*}%
The inequality $s_{(3,3)}\geq 0$ is equivalent to 
\begin{equation*}
-P^{2}QR+P^{2}Q+P^{2}R-2P+1\geq 0.
\end{equation*}

\subsection{Inequalities from monomial symmetric polynomials}

Let $m_{\lambda }$ denote the monomial symmetric polynomial associated with
the partition $\lambda $. We use the expansions 
\begin{align*}
m_{33}& =\sum_{\mathrm{sym}%
}x_{i}^{3}x_{j}^{3}=3e_{1}^{2}e_{4}-3e_{1}e_{2}e_{3}-3e_{1}e_{5}+e_{2}^{3}-3e_{2}e_{4}+3e_{3}^{2}+3e_{6},
\\
m_{222}& =\sum_{\mathrm{sym}%
}x_{i}^{2}x_{j}^{2}x_{k}^{2}=e_{3}^{2}-2e_{2}e_{4}+2e_{1}e_{5}-2e_{6}.
\end{align*}%
Hence 
\begin{equation*}
-12P^{2}QR+9P^{2}Q+6P^{2}R-6P+2=e_{2}^{-3}\left( 2m_{33}+3m_{222}\right)
\geq 0.
\end{equation*}

A further inequality, which is important later, follows from the monomial
expansion 
\begin{align}
& e_{2}e_{3}^{2}+e_{1}e_{3}e_{4}-2e_{4}e_{2}^{2}  \label{Gap} \\
&
=m_{332}+2m_{3221}+m_{32111}+4m_{2222}+6m_{22211}+6m_{221111}+5m_{2111111}%
\geq 0.  \notag
\end{align}%
Therefore, 
\begin{equation*}
e_{2}e_{3}^{2}(1+PR-2R)\geq 0,\qquad \text{or equivalently}\qquad 1\geq
(2-P)R.
\end{equation*}

\begin{rema}
The nonnegativity in \eqref{Gap0} and \eqref{Gap} does not follow directly
from the Newton--Maclaurin inequalities. The change-of-basis formulas above
allow us to use positivity in the power-sum and monomial bases.
\end{rema}

\subsection{Collection of inequalities}

The Newton inequalities also give 
\begin{equation*}
0<P\leq\frac{2(k-2)}{3(k-1)},\qquad 0<Q\leq\frac{k-1}{2k},\qquad 0<R\leq%
\frac{3(k-3)}{4(k-2)}.
\end{equation*}
We collect the above inequalities for later reference.

\begin{lemm}
\label{basis} Suppose that all the variables are positive. Then the
following inequalities hold:

\begin{enumerate}
\item \label{UpperQ} $1\geq(4-3P)Q$.

\item \label{UpperR} $1\geq(2-P)R$.

\item \label{LowerQR} $2P^2QR\geq2P-1$.

\item \label{Schur} $-P^2QR+P^2Q+P^2R-2P+1\geq0$.

\item \label{Mon} $-12P^2QR+9P^2Q+6P^2R-6P+2\geq0$.

\item \label{PQR} 
\begin{equation*}
0<P\leq\frac{2(k-2)}{3(k-1)},\qquad 0<Q\leq\frac{k-1}{2k},\qquad 0<R\leq%
\frac{3(k-3)}{4(k-2)}.
\end{equation*}
\end{enumerate}
\end{lemm}

\section{Reduction of the main estimate}

We first reduce the main estimate to a simpler algebraic inequality. Our
goal is to prove 
\begin{equation*}
M^{ij}\xi_i\xi_j =\frac{1}{F}\left(\sum_i F^{ii}\xi_i\right)^2
-\sum_{i,j}F^{ii,jj}\xi_i\xi_j +2\sum_{i>1}\frac{F^{ii}\xi_i^2}{\lambda_1}
-(1+\varepsilon)\frac{F^{11}\xi_1^2}{\lambda_1}\geq 0.
\end{equation*}

\subsection{Reduction to a determinant}

Since $M=(M^{ij})$ has at most one negative eigenvalue, it suffices to show
that $\det M\geq 0$. We compare $M$ with the simpler matrix 
\begin{equation*}
N=%
\begin{pmatrix}
A & B \\ 
C & D%
\end{pmatrix}%
,\qquad \det N=\det D\,(A-BD^{-1}C),
\end{equation*}%
where 
\begin{align*}
A& =-F^{11,11}+\frac{1}{F}(F^{11})^{2}-(1+\varepsilon )\frac{F^{11}}{\lambda
_{1}}, \\
B& =C^{T}=\left[ -F^{11,jj}+\frac{1}{F}F^{11}F^{jj}\right] _{2\leq j\leq n},
\\
D& =\func{diag}\left( \frac{2}{\lambda _{1}}\mathbf{v}\right) +\frac{1}{F}%
\mathbf{v}\mathbf{v}^{T}, \\
\mathbf{v}& =(F^{22},\ldots ,F^{nn})^{T}.
\end{align*}%
We note that $M$ and $N$ only differ at the lower-right block. 
\begin{equation*}
M=N+%
\begin{pmatrix}
0 & 0 \\ 
0 & D_2%
\end{pmatrix}%
,\qquad D_{2}=\left[ -F^{ii,jj}\right] _{2\leq i,j\leq n}\geq 0.
\end{equation*}%
Thus, it is enough to prove $N\geq 0.$ Since $D>0$, this is equivalent to 
\begin{equation*}
A-BD^{-1}C \geq 0.
\end{equation*}
To simplify the formulas, set 
\begin{align*}
a_j&=\frac{1}{\sigma_k^2} \left(\sigma_{k-1;1j}^2-\sigma_{k;1j}%
\sigma_{k-2;1j}\right), \\
b_j&=\frac{1}{\sigma_{k-1}^2} \left(\sigma_{k-2;1j}^2-\sigma_{k-1;1j}%
\sigma_{k-3;1j}\right), \\
m&=n-k.
\end{align*}

\subsection{Estimates for the matrix entries}

Throughout this section, the constants implicit in the $O(\cdot)$ notation may
depend on $n$, $k$, and $F$, but are independent of $\lambda_1$.

\paragraph{\textbf{The derivatives }$F^{ii}$.}

By Lemma \ref{Fijkl}, we have 
\begin{equation*}
F^{ii}=\frac{\sigma _{k-1;i}}{\sigma _{k-1}}-\frac{\sigma _{k}\sigma _{k-2;i}%
}{\sigma _{k-1}^{2}}=\frac{\sigma _{k-1;i}^{2}-\sigma _{k;i}\sigma _{k-2;i}}{%
\sigma _{k-1}^{2}}.
\end{equation*}%
Since $\lambda _{k}$ is bounded, the Newton inequality gives, for $j\geq 2$, 
\begin{align}
F^{jj}& =\frac{(\lambda _{1}\sigma _{k-2;1j}+\sigma _{k-1;1j})^{2}-(\lambda
_{1}\sigma _{k-1;1j}+\sigma _{k;1j})(\lambda _{1}\sigma _{k-3;1j}+\sigma
_{k-2;1j})}{\sigma _{k-1}^{2}}  \notag \\
& =\lambda _{1}^{2}b_{j}\left( 1+O\left( \frac{1}{\lambda _{1}}\right)
\right) .  \label{Fjj}
\end{align}%
Also, by Lemma \ref{Sigma_k-Lemma-0}, 
\begin{align}
\sum_{j=2}^{n}a_{j}& =\frac{1}{\sigma _{k}^{2}}\sum_{j=2}^{n}\left( \sigma
_{k-1;1j}^{2}-\sigma _{k;1j}\sigma _{k-2;1j}\right) =\frac{1}{\sigma _{k}^{2}%
}\sum_{j=2}^{n}\left( \sigma _{k-1;1}\sigma _{k-1;j}-\sigma _{k}\sigma
_{k-2;1j}\right)  \notag \\
& =\frac{1}{\sigma _{k}^{2}}\left( m\sigma _{k-1;1}^{2}-(m+1)\sigma
_{k;1}\sigma _{k-2;1}\right) .  \label{sum_aj}
\end{align}%
Therefore, 
\begin{equation}
F^{11}>\frac{F^{2}}{m}\sum_{j=2}^{n}a_{j}.  \label{F11(aj)}
\end{equation}%
Moreover, 
\begin{align}
F^{11}& =O\left( \frac{\left( \lambda _{2}...\lambda _{k}\right) ^{2}}{%
\left( \lambda _{1}...\lambda _{k-1}\right) ^{2}}\right) =O\left( \frac{1}{%
\lambda _{1}^{2}}\right) ,  \label{F11} \\
F^{nn}& =O\left( \frac{\left( \lambda _{1}...\lambda _{k-1}\right) ^{2}}{%
\left( \lambda _{1}...\lambda _{k-1}\right) ^{2}}\right) =O(1).  \label{Fnn}
\end{align}

\paragraph{\textbf{The term} $A$.}

We have 
\begin{align*}
F^{11,11}& =-2\frac{\sigma _{k-2;1}}{\sigma _{k-1}}F^{11}=\frac{-2\sigma
_{k-2;1}}{\lambda _{1}\sigma _{k-2;1}+\sigma _{k-1;1}}F^{11} \\
& =-2\frac{F^{11}}{\lambda _{1}}\left( 1+O\left( \frac{1}{\lambda _{1}}%
\right) \right) .
\end{align*}%
Thus, by \eqref{F11} and \eqref{F11(aj)}, 
\begin{align}
A& =-F^{11,11}+\frac{1}{F}\left( F^{11}\right) ^{2}-\left( 1+\varepsilon
\right) \frac{F^{11}}{\lambda _{1}}  \notag \\
& =(1-\varepsilon )\frac{F^{11}}{\lambda _{1}}+O\left( \frac{1}{\lambda
_{1}^{4}}\right)  \notag \\
& >(1-\varepsilon )\frac{F^{2}}{m\lambda _{1}}\sum_{j=2}^{n}a_{j}+O\left( 
\frac{1}{\lambda _{1}^{4}}\right) .  \label{A}
\end{align}

\paragraph{\textbf{The term} $B$.}

A direct computation gives 
\begin{equation}
B=-F^{11,jj}+\frac{1}{F}F^{11}F^{jj}=F(a_{j}-b_{j}).  \label{B}
\end{equation}

\paragraph{\textbf{The matrix} $D^{-1}$.}

Let 
\begin{equation*}
\Delta=\func{diag}\left(\frac{2}{\lambda_1}\mathbf{v}\right).
\end{equation*}
By the Sherman--Morrison formula, 
\begin{equation*}
D^{-1} =\left(\Delta+\frac{1}{F}\mathbf{v}\mathbf{v}^T\right)^{-1}
=\Delta^{-1} -\frac{\Delta^{-1}\mathbf{v}\mathbf{v}^T\Delta^{-1}} {F+\mathbf{%
v}^T\Delta^{-1}\mathbf{v}}.
\end{equation*}
Using \eqref{Fnn} and \eqref{Fjj}, we obtain 
\begin{align}
(D^{-1})_{ij} &=\frac{\lambda_1}{2F^{ii}}\delta_{ij} -\frac{\lambda_1^2}{%
4F+2\lambda_1\sum_{k=2}^nF^{kk}}  \notag \\
&=\frac{1}{2\lambda_1b_j}\delta_{ij} \left(1+O\left(\frac{1}{\lambda_1}%
\right)\right) -\frac{1}{2\lambda_1\sum_{k=2}^n b_k} \left(1+O\left(\frac{1}{%
\lambda_1}\right)\right).  \label{D-1}
\end{align}

\paragraph{\textbf{The term} $-BD^{-1}C$.}

By \eqref{B} and \eqref{D-1}, 
\begin{align}
-BD^{-1}C &=-\sum_{j=2}^n\frac{F^2(a_j-b_j)^2}{2\lambda_1b_j} +\frac{%
F^2\left(\sum_{j=2}^n(a_j-b_j)\right)^2} {2\lambda_1\sum_{j=2}^n b_j}
+O\left(\frac{1}{\lambda_1^4}\right)  \notag \\
&=\frac{F^2}{2\lambda_1} \left(-\sum_{j=2}^n\frac{a_j^2}{b_j} +\frac{%
\left(\sum_{j=2}^n a_j\right)^2} {\sum_{j=2}^n b_j}\right) +O\left(\frac{1}{%
\lambda_1^4}\right).  \label{BDC}
\end{align}

\subsection{The key inequality}

Combining \eqref{A} and \eqref{BDC}, we obtain 
\begin{align}
& A-BD^{-1}C  \notag \\
& >(1-\varepsilon )\frac{F^{2}}{m\lambda _{1}}\sum_{j=2}^{n}a_{j}\quad +%
\frac{F^{2}}{2\lambda _{1}}\left( -\sum_{j=2}^{n}\frac{a_{j}^{2}}{b_{j}}+%
\frac{\left( \sum_{j=2}^{n}a_{j}\right) ^{2}}{\sum_{j=2}^{n}b_{j}}\right)
+O\left( \frac{1}{\lambda _{1}^{4}}\right)  \notag \\
& =\frac{F^{2}}{2\lambda _{1}}\sum_{j=2}^{n}a_{j}\left( (1-\varepsilon )%
\frac{2}{m}-\frac{a_{j}}{b_{j}}+\frac{\sum_{k=2}^{n}a_{k}}{%
\sum_{k=2}^{n}b_{k}}\right) +O\left( \frac{1}{\lambda _{1}^{4}}\right) .
\label{Schurcomplement}
\end{align}%
Thus, the main estimate reduces to the following inequality.

\begin{lemm}
\label{lemma1} Suppose there exists a constant $\varepsilon _{0}(n,k)>0$
such that 
\begin{equation}
\frac{2}{m}-\frac{a_{j}}{b_{j}}+\frac{\sum_{k=2}^{n}a_{k}}{%
\sum_{k=2}^{n}b_{k}}>\varepsilon _{0}(n,k),\qquad 2\leq j\leq n.  \label{Key}
\end{equation}%
Then there exists $M=M(\varepsilon ,n,k,F)$ such that 
\begin{equation*}
\det N\geq 0
\end{equation*}%
whenever $\lambda _{1}>M$.
\end{lemm}

Indeed, choose $\varepsilon>0$ sufficiently small so that \eqref{Key}
implies
\[
(1-\varepsilon)\frac{2}{m}-\frac{a_j}{b_j}
+\frac{\sum_{k=2}^n a_k}{\sum_{k=2}^n b_k}
\geq \frac{\varepsilon_0}{2},
\qquad 2\leq j\leq n.
\]
Moreover,
\[
\sum_{j=2}^n a_j \geq a_n \geq \frac{c(n,k,F)}{\lambda_1^2}.
\]
It follows from \eqref{Schurcomplement} that
\[
A-BD^{-1}C
\geq \frac{c(n,k,F)}{\lambda_1^3}
-O\left(\frac{1}{\lambda_1^4}\right)>0
\]
when $\lambda_1$ is sufficiently large. Therefore $\det N\geq0$.

\begin{rema}
In fact, we only require $j=n$ in the inequality \eqref{Key}, since $%
a_{j}/b_{j}$ is increasing in $j$. Indeed, differentiating the quotient $%
\left( \sigma _{k-1}^{2}-\sigma _{k}\sigma _{k-2}\right) /\left( \sigma
_{k-2}^{2}-\sigma _{k-1}\sigma _{k-3}\right) $ with respect to $\lambda _{i}$
and using Lemma~\ref{basis}(\ref{Schur}) shows that the derivative is
positive. Also, we remark that \eqref{Key} can also be written as a weighted
variance.
\end{rema}

\section{Proof of the key inequality: endpoint cases}

It remains to prove the key inequality (\ref{Key}) used in the previous
section. We first introduce some notation that will be used throughout the
proof. Set 
\begin{equation}
e_{l}=\sigma _{l;1j},\qquad x=\lambda _{j}.  \label{notation}
\end{equation}%
Thus, $e_{l}$ is the $l$-th elementary symmetric polynomial in the remaining 
$n-2$ variables. We will also use 
\begin{equation}
\frac{\sigma _{k-1}^{2}}{\sigma _{k}^{2}}=\left( \frac{\sigma _{k-2;1}}{%
\sigma _{k-1;1}}\right) ^{2}\left( 1+O\left( \frac{1}{\lambda _{1}}\right)
\right) .  \label{simp}
\end{equation}%
Since all the terms involved are uniformly bounded, the error introduced by
(\ref{simp}) is $O(\lambda _1^{-1})$. Hence, for $\lambda _1$ sufficiently
large, it can be absorbed into the small constant $\varepsilon _0$ in
(\ref{Key}).

We first prove \eqref{Key} directly in the endpoint cases $k=2$ and
$k=n-1$. The remaining cases $3\leq k\leq n-2$ are treated in the
following sections.

\subsection{The case $k=2$}

Using (\ref{notation}) and (\ref{simp}), the left-hand side of (\ref{Key})
for $k=2$ becomes%
\begin{align}
&\frac{2}{m}-\frac{e_{1}^{2}-e_{2}}{\left( x+e_{1}\right) ^{2}}+\frac{1}{%
\left( x+e_{1}\right) ^{2}}\frac{m\left( x+e_{1}\right) ^{2}-\left(
m+1\right) \left( xe_{1}+e_{2}\right) }{\left( m+1\right) }  \notag \\
&=\frac{\left( m^{2}+2m+2\right) x^{2}+\left( m^{2}+3m+4\right)
e_{1}x+\left( m+2\right) e_{1}^{2}}{m\left( m+1\right) \left( x+e_{1}\right)
^{2}}.  \label{k=2}
\end{align}%
All the coefficients in the numerator are strictly positive. Hence, for
some $\varepsilon_0=\varepsilon_0(n,k)>0$, they remain strictly positive
after replacing $\frac{2}{m}$ by $\frac{2-\varepsilon_0}{m}$.
Therefore, \eqref{Key} holds for $k=2$.

\subsection{The case $k=n-1$}

It suffices to consider the first two terms in (\ref{Key}). We have%
\begin{align*}
&\frac{2}{m}-\left( \frac{\sigma _{k-2;1}}{\sigma _{k-1;1}}\right) ^{2}%
\frac{\sigma _{k-1;1j}^{2}-\sigma _{k;1j}\sigma _{k-2;1j}}{\sigma
_{k-2;1j}^{2}-\sigma _{k-1;1j}\sigma _{k-3;1j}} \\
&=2-\left( \frac{xe_{n-4}+e_{n-3}}{xe_{n-3}+e_{n-2}}\right) ^{2}\frac{%
e_{n-2}^{2}}{e_{n-3}^{2}-e_{n-2}e_{n-4}} \\
&=\frac{N\left( x\right) }{\left( xe_{n-3}+e_{n-2}\right)^2 \left(
e_{n-3}^{2}-e_{n-2}e_{n-4}\right) }\text{,}
\end{align*}%
where%
\begin{align*}
N(x) &=\left(
2e_{n-3}^{4}-2e_{n-3}^{2}e_{n-2}e_{n-4}-e_{n-2}^{2}e_{n-4}^{2}\right) x^{2}+
\\
&\quad \left( 4e_{n-3}^{3}e_{n-2}-6e_{n-3}e_{n-2}^{2}e_{n-4}\right) x+\left(
e_{n-3}^{2}e_{n-2}^{2}-2e_{n-2}^{3}e_{n-4}\right) .
\end{align*}%
By the Newton inequality in $n-2$ variables, we have%
\begin{equation*}
e_{n-3}^{2}\geq 2\frac{\left( n-2\right) }{\left( n-3\right) }e_{n-2}e_{n-4}%
\text{.}
\end{equation*}%
Therefore, all the coefficients of $N(x)$ are strictly positive. Hence,
for some $\varepsilon_0=\varepsilon_0(n,k)>0$, they remain strictly
positive after replacing $\frac{2}{m}$ by
$\frac{2-\varepsilon_0}{m}$. Thus, \eqref{Key} holds for $k=n-1$.

Therefore, (\ref{Key}) holds in both endpoint cases. It remains to consider $%
3\leq k\leq n-2$.

\section{Reduction to a cubic polynomial}

We now prove (\ref{Key}) for $3\leq k\leq n-2$. We first reduce the
inequality to the positivity of a cubic polynomial. By (\ref{sum_aj}), the
left-hand side of (\ref{Key}) can be written as 
\begin{equation}
\frac{2}{m}-\frac{\sigma _{k-1}^{2}}{\sigma _{k}^{2}}\left( \frac{\sigma
_{k-1;1j}^{2}-\sigma _{k;1j}\sigma _{k-2;1j}}{\sigma _{k-2;1j}^{2}-\sigma
_{k-1;1j}\sigma _{k-3;1j}}-\frac{m\sigma _{k-1;1}^{2}-(m+1)\sigma
_{k;1}\sigma _{k-2;1}}{(m+1)\sigma _{k-2;1}^{2}-(m+2)\sigma _{k-1;1}\sigma
_{k-3;1}}\right) .  \label{5.1}
\end{equation}%
Using (\ref{notation}) and (\ref{simp}), we write (\ref{5.1}) in terms of $x$
and $e_{l}$. After taking a common denominator and canceling common factors,
its numerator is a cubic polynomial in $x$: 
\begin{equation}
N(x)=A_{3}x^{3}+A_{2}x^{2}+A_{1}x+A_{0}.  \label{cubic}
\end{equation}%
Set 
\begin{equation*}
\mathcal{P}=\frac{e_{k-3}e_{k-1}}{e_{k-2}^{2}},\qquad \mathcal{Q}=\frac{%
e_{k-2}e_{k}}{e_{k-1}^{2}},\qquad \mathcal{R}=\frac{e_{k-2}e_{k-4}}{%
e_{k-3}^{2}}.
\end{equation*}%
Then the coefficients in (\ref{cubic}) are given by 
\begin{align*}
A_{0}& =\frac{\mathcal{P}e_{k-2}^{6}}{e_{k-3}}\Big(2(m+2)\mathcal{P}^{2}-(m%
\mathcal{Q}+2m+6)\mathcal{P}+(m+2)\Big), \\
A_{1}& =e_{k-2}^{5}\Big(2(m+2)\mathcal{R}\mathcal{P}^{3}+\big(m^{2}(\mathcal{%
Q}+\mathcal{R}-\mathcal{Q}\mathcal{R})-m(2\mathcal{Q}\mathcal{R}+\mathcal{Q}%
-4)+4(1-\mathcal{R})\big)\mathcal{P}^{2} \\
& \qquad -(2m^{2}+5m+6)\mathcal{P}+(m^{2}+2m+2)\Big), \\
A_{2}& =e_{k-3}e_{k-2}^{4}\Big(\big(2m^{2}(\mathcal{Q}+\mathcal{R}-\mathcal{Q%
}\mathcal{R})+m(\mathcal{Q}+8\mathcal{R}-4\mathcal{Q}\mathcal{R})+8\mathcal{R%
}-2\big)\mathcal{P}^{2} \\
& \qquad -\big(4m^{2}+3m-2+4m\mathcal{R}+8\mathcal{R}\big)\mathcal{P}+2m(m+1)%
\Big), \\
A_{3}& =e_{k-3}^{2}e_{k-2}^{3}\Big(\big(m^{2}(\mathcal{Q}+\mathcal{R}-%
\mathcal{Q}\mathcal{R})+m(\mathcal{Q}+2\mathcal{R}-2\mathcal{Q}\mathcal{R})%
\big)\mathcal{P}^{2} \\
& \qquad +\big(2m\mathcal{R}+4\mathcal{R}-2m^{2}-3m-2\big)\mathcal{P}+\big(%
m^{2}+2m+2-2m\mathcal{R}-4\mathcal{R}\big)\Big).
\end{align*}

We next rewrite these coefficients in terms of the variables $P,Q,R$
introduced in Section~3. Let $E_j$ denote the normalized elementary symmetric
polynomial. We use the following two standard identities. If
$X=(x_1,\ldots,x_n)$, then:

\begin{itemize}
\item By Rolle's theorem, there exists $Y=(y_1,\ldots,y_{n-1})$ such that 
\begin{equation}
E_j(X)=E_j(Y), \qquad 0\leq j\leq |Y|.  \label{rolle's}
\end{equation}

\item If $X^{-1}=\{1/a:a\in X\}$, then 
\begin{equation}
E_j(X^{-1})=\frac{E_{n-j}(X)}{E_n(X)}, \qquad 0\leq j\leq n.  \label{recipr}
\end{equation}
\end{itemize}

Applying (\ref{rolle's}) repeatedly, we reduce from $n-2$ to $k$ variables
while preserving the normalized elementary symmetric polynomials $E_j$,
$0\leq j\leq k$. We then pass to reciprocal variables using (\ref{recipr}),
which converts the high-degree quantities $E_{k-j}$ into the low-degree
quantities $E_j$.

Let 
\begin{equation*}
\lambda ^{\prime }=(\lambda _{1}^{\prime },\ldots ,\lambda _{n-2}^{\prime }).
\end{equation*}%
Then there exists $\mu =(\mu _{1},\ldots ,\mu _{k})$ such that 
\begin{align*}
\mathcal{P} &=\frac{e_{k-3}e_{k-1}}{e_{k-2}^{2}}=\frac{m}{m+1}\frac{k-2}{k-1%
}\frac{E_{k-3}(\lambda ^{\prime })E_{k-1}(\lambda ^{\prime })}{%
E_{k-2}^{2}(\lambda ^{\prime })}=\frac{m}{m+1}\frac{k-2}{k-1}\frac{%
E_{k-3}(\mu )E_{k-1}(\mu )}{E_{k-2}^{2}(\mu )} \\
&=\frac{m}{m+1}\frac{k-2}{k-1}\frac{E_{3}(\mu ^{-1})E_{1}(\mu ^{-1})}{%
E_{2}^{2}(\mu ^{-1})}=\frac{3m}{2(m+1)}\frac{e_{1}e_{3}}{e_{2}^{2}}.
\end{align*}%
Similarly, 
\begin{equation*}
\mathcal{Q}=\frac{2(m-1)}{m}\frac{e_{2}}{e_{1}^{2}},\qquad \mathcal{R}=\frac{%
4(m+1)}{3(m+2)}\frac{e_{2}e_{4}}{e_{3}^{2}}.
\end{equation*}%
We still denote the elementary symmetric polynomials
of $\mu^{-1}$ by $e_j$.
Recall that 
\begin{equation*}
P=\frac{e_{1}e_{3}}{e_{2}^{2}},\qquad Q=\frac{e_{2}}{e_{1}^{2}},\qquad R=%
\frac{e_{2}e_{4}}{e_{3}^{2}}.
\end{equation*}%
Therefore, 
\begin{equation}
\mathcal{P}=\frac{3m}{2(m+1)}P,\qquad \mathcal{Q}=\frac{2(m-1)}{m}Q,\qquad 
\mathcal{R}=\frac{4(m+1)}{3(m+2)}R.  \label{PQRsub}
\end{equation}
Thus the coefficients $A_i$ can be expressed entirely in terms of $P,Q,R$
and $m$.

Substituting these expressions into $A_i$, we obtain four polynomials in $%
e_1,e_2,e_3,e_4$. For example, up to a positive factor, 
\begin{equation*}
A_0= 9m^{2}(m+2)e_1^{3}e_3^{2} -6m(m+1)(m+3)e_1^{2}e_2^{2}e_3
-6m(m-1)(m+1)e_2^{3}e_3 +2(m+1)^{2}(m+2)e_1e_2^{4}.
\end{equation*}

The Newton--Maclaurin inequalities alone do not seem sufficient to prove the
positivity of $A_0$. The dominant term, $e_1e_2^4$, has a coefficient that
is too small. Another positive term, $e_1^3e_3^2$, is dominated by the
negative term $e_1^2e_2^2e_3$, while its comparison with $e_2^3e_3$ is
unclear---it can be either larger or smaller. This suggests that the
elementary symmetric polynomials may not be a suitable basis.

We first give an independent computational verification of the positivity
of $A_i$. A complete direct proof using the symmetric polynomial inequalities
established earlier is given in Section~8.

\section{A computational verification}

We give an independent computational verification of the positivity of
$A_i$. Using SageMath \cite{S}, we expand $A_i$ in the
monomial basis. All coefficients in the expansions of $A_2$ and $A_3$
are positive. For $A_0$ and $A_1$, only a few coefficients are negative,
and they can be controlled by positive terms using Muirhead's inequality.
For example, the only negative term in the expansion of $A_1$ is
$m_{2^{2}1^{6}}$, which is controlled by $m_{31^{7}}$ since
\[
7!m_{31^{7}}\geq 2!6!m_{2^{2}1^{6}}.
\]
This verifies the positivity of $A_i$. A complete direct proof is given in
the next section.

\section{A direct proof of positivity}

We now give a direct proof of the positivity of the coefficients in (\ref%
{cubic}). Instead of working in the elementary symmetric basis, we use the
variables $P,Q,R$ and the symmetric polynomial inequalities established in
Section~3.

\subsection{Main idea}

The main observation is the \emph{multi-affine }structure. For fixed $%
\mathcal{P}$, each coefficient $A_{i}$ in (\ref{cubic}) is multi-affine in $%
\mathcal{Q}$ and $\mathcal{R}$. Therefore, if $(\mathcal{Q},\mathcal{R})$ is
restricted to a suitable rectangle, the minimum of $A_{i}$ is attained at
one of its vertices.

The case of $A_0$ is direct. For $A_1,A_2,A_3$, finding a suitable rectangle
is more difficult. We therefore write each $A_i$ as a polynomial in $m$ and
use the symmetric polynomial inequalities from Section~3 to obtain suitable
bounds for each coefficient.

Since only the sign of $A_i$ matters, from now on we allow $A_i$ to differ
from its previous definition by a positive factor.

\subsection{The coefficient $A_0$}

In terms of $P,Q,R$, and up to a positive factor, $A_{0}$ can be written as%
\begin{equation*}
A_{0}=e_{1}^{9}Q^{4}\,F_{0}(P,Q,R;m),
\end{equation*}%
where 
\begin{equation*}
F_{0}=9m^{2}(m+2)P^{2}-6m(m+1)(m+3)P-6m(m-1)(m+1)PQ+2(m+1)^{2}(m+2).
\end{equation*}%
Since $F_{0}$ is decreasing in $Q$, Lemma \ref{basis}(\ref{UpperQ}) shows
that its minimum occurs at $Q=1/(4-3P)$. Substituting this value, we obtain%
\begin{equation*}
F_{0}=\frac{\left( 2m+2-3Pm\right) }{4-3P}\left( 9m\left( m+2\right)
P^{2}-6\left( 2m^{2}+5m+1\right) P+4\left( m+1\right) \left( m+2\right)
\right) \text{.}
\end{equation*}%
The discriminant of the quadratic polynomial in the second bracket is 
\begin{equation*}
-36\left( 3m^{2}+6m-1\right) .
\end{equation*}%
Therefore
\[
F_0>C(n,k)>0
\qquad\text{for}\qquad
0\leq P\leq \frac{2(k-2)}{3(k-1)}.
\]
Since $F_0$ is a polynomial and the admissible set is contained in a
compact subset of
$[0,\frac{2}{3}]\times[0,\frac{1}{2}]\times[0,\frac{3}{4}]$,
this strict positivity is preserved under a sufficiently small perturbation
of the coefficient $\frac{2}{m}$. Hence, for some
$\varepsilon_0=\varepsilon_0(n,k)>0$, $F_0$ remains positive when
$\frac{2}{m}$ is replaced by $\frac{2-\varepsilon_0}{m}$.

\subsection{The coefficients $A_1,A_2,A_3$}

We write the three coefficients in the following form.%
\begin{equation*}
A_{1}=e_{1}^{10}Q^{5}\,F_{1}(P,Q,R;m),\text{ }A_{2}=e_{1}^{11}Q^{6}P%
\,F_{2}(P,Q,R;m)\text{, }A_{3}=e_{1}^{12}Q^{7}P^{2}\,F_{3}(P,Q,R;m)\text{,}
\end{equation*}%
where%
\begin{align*}
F_{1}& =\alpha m^{4}+\beta _{1}m^{3}+\gamma _{1}m^{2}+\left( 12-18P\right)
m+4. \\
F_{2}&=2\alpha m^{3}+\beta _{2}m^{2}+\gamma _{2}m-16PR+6P+4, \\
F_{3}& =3\alpha m^{3}+\beta _{3}m^{2}+\gamma _{3}m-16R+12.
\end{align*}%
By Lemma \ref{basis}(\ref{UpperR},\ref{PQR}), the sum of all terms in $F_{i}$
other than $\alpha $, $\beta _{i}$, and $\gamma _{i}\,\ $has a positive
lower bound. Thus, it is enough to prove that $\alpha $, $\beta _{i}$, and $%
\gamma _{i}$ are nonnegative.

\subsubsection{The term $\protect\alpha$}

\begin{equation*}
\alpha =-12P^{2}QR+9P^{2}Q+6P^{2}R-6P+2,
\end{equation*}%
which is nonnegative by Lemma \ref{basis}(\ref{Mon}).

\subsubsection{The terms $\protect\beta_i$}

We will show that 
\begin{align*}
\beta _{1} &=18P^{3}R-18P^{2}Q-6P^{2}R+18P^{2}-21P+8\geq 0, \\
\beta _{2} &=-9P^{2}Q+36P^{2}R-16PR-21P+12\geq 0, \\
\beta _{3} &=36P^{2}QR-27P^{2}Q+24PR-27P-16R+18\geq 0.
\end{align*}%
We first establish a lower bound for $R$. Since $\alpha \geq0$, by Lemma \ref%
{basis}(\ref{PQR}), 
\begin{equation*}
\partial _{Q}\alpha =3P^{2}(3-4R)>0.
\end{equation*}%
By Lemma \ref{basis}(\ref{UpperQ}), setting $Q=1/(4-3P)$ in $\alpha $ gives 
\begin{equation*}
\alpha \left( P,\frac{1}{4-3P},R\right) =\frac{2-3P}{4-3P}\left(
6P^{2}R-9P+4\right) \geq 0.
\end{equation*}%
Therefore, 
\begin{equation*}
R\geq\frac{9P-4}{6P^{2}}.
\end{equation*}%
Combining this with Lemma \ref{basis}(\ref{UpperR}), we obtain 
\begin{equation*}
\frac{1}{2-P}\geq R\geq\frac{9P-4}{6P^{2}}.
\end{equation*}%
Next we use the above bounds to place the admissible region inside a rectangle. Since the functions are multi-affine in $Q$ and $R$, their minimum on the rectangle is attained at a vertex. We consider the three terms separately.
\medskip

\paragraph{\textbf{For }$\protect\beta _{1}$}

Observe that 
\begin{equation*}
\partial _{Q}\beta _{1}=-18P^{2}<0,\qquad \partial _{R}\beta
_{1}=6P^{2}(3P-1).
\end{equation*}%
The sign of $\partial_R\beta_1$ changes at $P=1/3$. 
\begin{equation*}
\begin{array}{ccc}
& \text{Minimum vertex }(Q,R) & \text{Minimum value of }\beta _{1} \\%
[2mm] \hline
&  &  \\[-2mm] 
0<P\leq \dfrac{1}{3} & \left( \dfrac{1}{4-3P},\dfrac{1}{2-P}\right) & \dfrac{%
(2-3P)(4-5P)(8-9P)}{(2-P)(4-3P)} \\[4mm] 
\dfrac{1}{3}<P\leq \dfrac{2\left( k-2\right) }{3(k-1)} & \left( \dfrac{1}{%
4-3P},\dfrac{9P-4}{6P^{2}}\right) & \dfrac{3(2-3P)^{2}(4-5P)}{4-3P}%
\end{array}%
\end{equation*}%
The minimum of $\beta _{1}$ on each interval is positive by direct calculus.
\medskip
\paragraph{\textbf{For }$\protect\beta _{2}$}

\begin{equation*}
\partial _{Q}\beta _{2}=-9P^{2}<0,\qquad \partial _{R}\beta _{2}=4P(9P-4).
\end{equation*}%
$\partial _{R}\beta _{2}$ changes sign at $P=4/9$.
\begin{equation*}
\begin{array}{ccc}
& \text{Minimum vertex }(Q,R) & \text{Minimum value of }\beta _{2} \\%
[2mm] \hline
&  &  \\[-2mm] 
0<P\leq \dfrac{4}{9} & \left( \dfrac{1}{4-3P},\dfrac{1}{2-P}\right) & \dfrac{%
2(2-3P)(27P^{2}-52P+24)}{(2-P)(4-3P)} \\[4mm] 
\dfrac{4}{9}<P\leq \dfrac{2\left( k-2\right) }{3(k-1)} & \left( \dfrac{1}{%
4-3P},\dfrac{9P-4}{6P^{2}}\right) & \dfrac{4(3P-2)^{2}(9P-8)}{3P(3P-4)}%
\end{array}%
\end{equation*}%
The minimum of $\beta _{2}$ on each interval is positive by direct calculus.
\medskip
\paragraph{\textbf{For }$\protect\beta _{3}$}

\begin{equation*}
\partial _{Q}\beta _{3}=9P^{2}(4R-3)<0,\qquad \left. \partial _{R}\beta
_{3}\right\vert _{Q=\frac{1}{4-3P}}=\frac{4(9P^{2}-36P+16)}{3P-4}
\end{equation*}%
$\partial _{R}\beta _{3}$ changes sign at $P_{0}=2-\frac{2}{3}\sqrt{5}$. 
\begin{equation*}
\begin{array}{ccc}
& \text{Minimum vertex }(Q,R) & \text{Minimum value of }\beta _{3} \\%
[2mm] \hline
&  &  \\[-2mm] 
0<P\leq P_{0} & \left( \dfrac{1}{4-3P},\dfrac{1}{2-P}\right) & \dfrac{%
2(2-3P)(9P^{2}-33P+20)}{(2-P)(4-3P)} \\[4mm] 
P_{0}<P\leq \dfrac{2\left( k-2\right) }{3(k-1)} & \left( \dfrac{1}{4-3P},%
\dfrac{9P-4}{6P^{2}}\right) & \dfrac{2(4-3P)(2-3P)^{2}}{3P^{2}}%
\end{array}%
\end{equation*}%
The minimum of $\beta _{3}$ on each interval is positive by direct calculus.

\subsubsection{The terms $\protect\gamma_i$}

We will show that 
\begin{align*}
\gamma _{1} &=12P^{2}QR+9P^{2}Q-12P^{2}R+18P^{2}-33P+14\geq 0, \\
\gamma _{2} &=24P^{2}QR-9P^{2}Q+24P^{2}R-9P^{2}-32PR-3P+12\geq 0, \\
\gamma _{3} &=24PR-18P-32R+24\geq 0.
\end{align*}%
We first establish a lower bound for $Q$. From Lemma \ref{basis}(\ref{UpperR}%
,\ref{LowerQR}),%
\begin{equation*}
0\leq 1-2P+2P^{2}QR\leq 1-2P+\frac{2P^{2}Q}{2-P}\text{,}
\end{equation*}%
we obtain 
\begin{equation*}
Q\geq Q_{0}:=\frac{(2P-1)(2-P)}{2P^{2}}.
\end{equation*}%
Combining this with Lemma \ref{basis}(\ref{UpperQ}), we obtain 
\begin{equation*}
\dfrac{1}{4-3P}\geq Q\geq \max \left\{ Q_{0},0\right\} .
\end{equation*}%
We first consider $\gamma _{1}$.%
\begin{equation*}
\partial _{R}\gamma _{1}=12P^{2}(Q-1)<0,\qquad \left. \partial _{Q}\gamma
_{1}\right\vert _{R=\frac{1}{2-P}}=\frac{3P^{2}(10-3P)}{2-P}>0
\end{equation*}%
Therefore, 
\begin{equation*}
\gamma _{1}\geq 
\begin{cases}
\gamma _{1}\!\left( P,0,\frac{1}{2-P}\right) =\frac{(1-2P)(9P^{2}-24P+28)}{%
2-P}, & 0<P\leq \frac{1}{2}, \\[2mm] 
\gamma _{1}\!\left( P,Q_{0},\frac{1}{2-P}\right) =\frac{(2P-1)(2-3P)(3P+2)}{%
4-2P}, & \frac{1}{2}<P<\frac{2}{3}.%
\end{cases}%
\end{equation*}%
Hence $\gamma _{1}\geq 0$.

\medskip

Next, consider $\gamma _{2}$. By Lemma \ref{basis}(\ref{PQR}),%
\begin{equation*}
\begin{cases}
\partial _{R}\gamma _{2}=24P^{2}(Q+1)-32P<36P^{2}-32P=4P(9P-8)<0, \\[2mm] 
\left. \partial _{Q}\gamma _{2}\right\vert _{R=\frac{1}{2-P}}=\frac{%
3P^{2}(3P+2)}{2-P}>0%
\end{cases}%
\end{equation*}%
Therefore,%
\begin{equation*}
\gamma _{2}\geq \gamma _{2}\!\left( P,Q_{0},\frac{1}{2-P}\right) =\frac{%
(2-3P)(18-17P)}{2(2-P)},\qquad 0<P<\frac{2}{3}
\end{equation*}%
Hence $\gamma _{2}>0$.

\medskip

Finally, consider $\gamma _3$. Since 
\begin{equation*}
\partial _{R}\gamma _{3}=24P-32<0,
\end{equation*}%
we have 
\begin{equation*}
\gamma _{3}\geq \gamma _{3}\!\left( P,\frac{1}{2-P}\right) =\frac{%
2(4-3P)(2-3P)}{2-P}>0.
\end{equation*}

We have proved that $\alpha$, $\beta_i$, and $\gamma_i$ are all
nonnegative. Since the remaining terms in $F_i$ have a positive lower
bound depending only on $n$ and $k$, we obtain
\[
F_i>C(n,k)>0,\qquad 1\leq i\leq3.
\]
Since each $F_i$ is a polynomial, the same compactness argument as above
shows that this positivity is preserved under a sufficiently small
perturbation of $\frac{2}{m}$. After decreasing
$\varepsilon_0(n,k)>0$ if necessary, we may therefore replace
$\frac{2}{m}$ by $\frac{2-\varepsilon_0}{m}$.
This completes the direct proof.

\section*{Acknowledgments}

The author would like to thank Professor Siyuan Lu for helpful discussions.


\begin{thebibliography}{99}

\bibitem{Chu}
J. Chu,
A simple proof of curvature estimate for convex solution of $k$-Hessian equation,
\emph{Proc. Amer. Math. Soc.} \textbf{149} (2021), no.~8, 3541--3552.

\bibitem{GRW}
P. Guan, C. Ren, and Z. Wang,
Global $C^2$ estimates for convex solutions of curvature equations,
\emph{Comm. Pure Appl. Math.} \textbf{68} (2015), no.~8, 1287--1325.

\bibitem{GS}
P. Guan and M. Sroka,
A special concavity property for positive Hessian quotient operators,
\emph{Discrete Contin. Dyn. Syst.} \textbf{54} (2026), 50--60.

\bibitem{LiWu}
Z. Li and K. Wu,
A concavity inequality and interior $C^2$ estimate for Hessian quotient equations,
preprint, arXiv:2608.17405.

\bibitem{LT}
M. Lin and N. S. Trudinger,
The Dirichlet problem for the prescribed curvature quotient equations,
\emph{Topol. Methods Nonlinear Anal.} \textbf{3} (1994), 307--323.

\bibitem{Lu24}
S. Lu,
Interior $C^2$ estimate for Hessian quotient equation in general dimension,
\emph{Ann. PDE} \textbf{11} (2025), no.~2, Paper No.~17.

\bibitem{LuHyperbolic}
S. Lu,
Curvature estimates for semi-convex solutions of Hessian equations in hyperbolic space,
\emph{Calc. Var. Partial Differential Equations} \textbf{62} (2023),
no.~9, Paper No.~257.

\bibitem{LuTsai}
S. Lu and Y.-L. Tsai,
A note on interior $C^2$ estimate for general Hessian quotient equation,
\emph{Commun. Pure Appl. Anal.} \textbf{33} (2026), 88--100.

\bibitem{M}
I.~G. Macdonald,
\emph{Symmetric Functions and Hall Polynomials},
2nd ed., Oxford Mathematical Monographs,
Clarendon Press, Oxford University Press, Oxford, 1995.

\bibitem{RW}
C. Ren and Z. Wang,
On the curvature estimates for Hessian equations,
\emph{Amer. J. Math.} \textbf{141} (2019), no.~5, 1281--1315.

\bibitem{RW2}
C. Ren and Z. Wang,
The global curvature estimate for the $n-2$ Hessian equation,
\emph{Calc. Var. Partial Differential Equations} \textbf{62} (2023),
no.~9, Paper No.~239.

\bibitem{RWNotes}
C. Ren and Z. Wang,
Notes on the curvature estimates for Hessian equations,
preprint, arXiv:2003.14234.

\bibitem{S}
The Sage Developers,
\emph{SageMath, the Sage Mathematics Software System},
\url{https://www.sagemath.org}.

\bibitem{Zhang}
R. Zhang,
$C^2$ estimates for $k$-Hessian equations and a rigidity theorem,
\emph{Adv. Math.} \textbf{480} (2025), Paper No.~110488.

\end{thebibliography}
\end{document}